\documentclass{ifacconf}

\usepackage{amsmath}
\usepackage{amssymb}
\usepackage{amsfonts}
\usepackage{physics}
\usepackage{enumerate}
\usepackage{xcolor}
\usepackage{balance}

\def\RR{\mathbb{R}}

\def\EE{\mathbb{E}}
\def\PP{\mathbb{P}}
\def\I{\mathcal{I}}

\def\U{\mathcal{U}}
\def\xhat{\hat{x}}

\def\bs{\boldsymbol}
\def\bxhat{\hat{\bf x}}
\def\bx{{\bf x}}
\def\bu{{\bf u}}
\def\bv{{\bf v}}
\def\bw{{\bf w}}
\def\Q{{\bf Q}}
\def\R{{\bf R}}
\def\ones{{\bf 1}}
\def\defeq{:=}
\def\tail{\circ}
\DeclareMathOperator{\diag}{diag}
\DeclareMathOperator{\vvec}{vec}
\newtheorem{assumption}{Assumption}

\def\bpf{\textnormal{\textit{Proof:}\hspace{1ex}}}
\def\epf{\hfill \mbox{\qed}}%$\Box$}
\usepackage{graphicx}      % include this line if your document contains figures
\usepackage{natbib}        % required for bibliography

\begin{document}
\begin{frontmatter}

\title{Output feedback stochastic MPC with packet losses}

% \thanks[footnoteinfo]{Sponsor and financial support acknowledgment goes here. Paper titles should be written in uppercase and lowercase letters, not all uppercase.}

\author[First]{Shuhao Yan}
\quad
\author{Mark Cannon}
\quad
\author{Paul Goulart}
% \author[Second]{Second B. Author, Jr.}
% \author[Third]{Third C. Author}

\address[First]{Department of Engineering Science, University of Oxford, UK\\ Email: \{shuhao.yan, mark.cannon, paul.goulart\}@eng.ox.ac.uk}
% \address[Second]{mark.cannon@eng.ox.ac.uk}
% \address[Third]{paul.goulart@eng.ox.ac.uk}

\begin{abstract}                % Abstract of not more than 250 words.
The paper considers constrained linear systems with stochastic additive disturbances and noisy measurements transmitted over a lossy communication channel. We propose a model predictive control (MPC) law that minimizes a discounted cost subject to a discounted expectation constraint.
%
%The cost is an infinite sum of discounted probabilities of violation of linear bounds on future states.
%, and both the cost and constraint are evaluated over an infinite horizon.
%v
Sensor data is assumed to be lost with known probability, and data losses are accounted for by expressing the predicted control policy as an affine function of future observations, which results in a convex optimal control problem.
An online constraint-tightening technique ensures recursive feasibility of the online optimization and satisfaction of the expectation constraint without bounds on the distributions of the noise and disturbance inputs.
The cost evaluated along trajectories of the closed loop system is shown to be bounded by the optimal predicted cost.
A numerical example is given to illustrate these results.
\end{abstract}


\end{frontmatter}
%===============================================================================

\section{Introduction}
Robust model predictive control often considers worst-case disturbance
bounds, so that hard constraints on system states and control inputs are satisfied for all possible disturbances \citep{mayne00,Mesbah16_review,kouvaritakis2015}. However, worst-case disturbance bounds can be extremely conservative or even non-existent, which motivates the development of stochastic MPC with chance constraints. In many applications of practical interest, system states cannot be measured directly and instead have to be estimated from output measurements. Existing stochastic MPC algorithms incorporating state estimation~\citep[e.g.][]{CANNON2012,Dai15} typically do not consider optimizing state feedback gains online, and the estimator gain is typically chosen as the steady state Kalman filter gain.

Control systems that rely on sensor signals transmitted over a network must tolerate communication delays and data losses. These pose additional challenges for estimation and control problems when constraints are present. From a control perspective, these features can be modelled as information losses by random processes, such as Bernoulli processes \citep{sinopoli04} or Markov chains \citep{Leong17CDC}.
In \citet{sinopoli04}, the arrival of output observations is modelled as a Bernoulli process and fundamental results are derived, including bounds on the critical value for the arrival probability of the observation update and convergence properties of the algebraic Riccati equation for Kalman filters with intermittent observations. In \citet{schenato07}, it is shown that the well-known separation principle holds with sensor packet losses, whereas this is not the case if constraints are present. \citet{mishra2019stochastic} consider the problem of controlling linear systems with unbounded additive disturbances and measurement noise by using an affine policy, where both sensor measurements and control actions are lost with given probabilities. 
Alternatively, these problems can be modelled as jump linear systems \citep{mariton1990} switching between different states according to a transition probability matrix.

This paper designs an output-feedback MPC algorithm to minimize a discounted cost function subject to a discounted expectation constraint, assuming sensor measurements to be lost with a given probability. The discount setting is common to many control problems \citep[e.g.][]{bertsekas1995dynamic,VanParys13,KOUVARITAKIS03,KAMGARPOUR17}, and an appropriate discounting factor can provide stability guarantees \citep{postoyan17}.
%An appropriate discounting factor can be helpful in deriving theoretical properties, while in some cases it needs to be tuned close to 1 to ensure stability guarantees \citep{postoyan17}.
In this work, the discounting factor allows consideration of unbounded disturbances and measurement noise,
%close to 1 to ensure stability guarantees \citep{postoyan17}. By the discounting factor,
and we derive bounds on the cost and constraints for the closed loop system using a constraint-tightening technique \citep{yan18}.
%In the controller design,
Instead of choosing the future control policy as pre-stabilising feedback with perturbations \citep{Cannon11}, we parameterise predicted control inputs as affine functions of future output measurements and show that the problem of optimizing the associated feedback gains is convex.
This allows the distributions of future states to be controlled even when output measurements are lost.
%Since in this setting, the former parameterisation can be very suboptimal, although it provides a more computationally efficient algorithm.

This paper is organised as follows. We describe the control problem in Section \ref{sec:problem description}, and introduce the controller parameterization and implementation in Section \ref{sec:control_parameterization}.
We compute predicted state and control sequences via their first and second moments in Section \ref{sec:MPC_predictions}. In Section \ref{sec:cost and constraints}, we derive the terminal conditions and give explicit expressions for the cost and constraints. Our main results, including a closed loop cost bound and constraint satisfaction, are in Section \ref{sec:properties}. Section~\ref{sec:example} provides a numerical example and the paper is concluded in Section \ref{sec:conclusion}.

\textit{Notation}:
The $n\times n$ identity matrix is $I_{n\times n}$, and the $n\times m$ matrix with all elements equal to $1$ is $\ones_{n\times m}$. 
%and ${\ones_n\defeq\ones_{n\times 1}}$.
%
The vectorized form of a matrix $A= [a_1 \ \cdots \ a_n]$ is $\vvec(A)\defeq[a_1^\top \ \cdots \ a_n^\top]^\top$ and $A\otimes B$ is the Kronecker product.
The Euclidean norm is $\| x \|$ and, for a matrix $Q$, $Q\succ 0$ ($Q\succeq 0$) indicates that $Q$ is positive definite (semidefinite) and $\|x\|_Q^2 \defeq x^\top Q x$.

\section{Problem description}\label{sec:problem description}

\subsection{System model and feedback information}

We assume a system with linear discrete time dynamics
\vspace{-1mm}
\begin{subequations}\label{eq:system}
\begin{align}
x_{k+1} &= A x_k + B u_k + D w_k,\\
y_k &= Cx_k + v_k,~
z_k = \gamma_k y_k
\end{align}
\end{subequations}
where $x\in\RR^{n_x}$, $u\in\RR^{n_u}$, $y\in\RR^{n_y}$, $z\in\RR^{n_y}$ are the state, control input, sensor measurement, and the measurement information received by the controller respectively.
The disturbance, measurement noise and packet loss sequences, $\{w_k\}_{k=0}^\infty$, $\{v_k\}_{k=0}^\infty$ and $\{\gamma_k\}_{k=0}^\infty$, are assumed to have independent, identically distributed (i.i.d.) elements with
\begingroup
\setlength{\abovedisplayskip}{-1pt}
\setlength{\belowdisplayskip}{-0.5pt}
\begin{alignat*}{2}
& \EE\{w_k\} = 0 , & \qquad & \EE\{w_kw_k^\top\} = \Sigma_w \succeq 0 , \\
& \EE\{v_k\} = 0 , & & \EE\{v_kv_k^\top\} = \Sigma_v \succ 0 ,\\
& \PP\{\gamma_k= 0\} = 1-\lambda ,  & & \PP\{\gamma_k = 1\} =  \lambda .
\end{alignat*}
\endgroup
The variable $\gamma_k\in\{0,1\}$ indicates whether sensor data at the $k$th sampling instant is received by the controller.
The information available to the controller at time $k$ consists of
$\{u_i\}_{i=0}^{k-1}$, $\{(z_i,\gamma_i)\}_{i=0}^k$,
the initial mean $\EE\{x_0\}=\xhat_0$, and covariance $\EE\{(x_0-\xhat_0)(x_0 - \xhat_0)^\top\} = \Sigma_0$ of the model state.

We define the information sets
\[
\I_{k} \defeq \{\I_{k-1}, (z_k,\gamma_k)\},
\quad
\U_{k} \defeq  \{\U_{k-1},u_k\},
\]
for all $k\geq 0$, where $\I_{-1} \defeq  \{\xhat_0,\Sigma_0\}$, $\U_{-1} \defeq \{\,\}$.  Finally, we define conditional expectation operators as
\[
\EE_k\{ \cdot \} \defeq \EE\{ \cdot \, | \, \U_{k-1},\I_{k-1}\},
\quad
\EE\{ \cdot \} \defeq \EE_0\{\cdot\}.
\]

\begin{assumption}\label{assumption:ctrb_and_obsv}
The pair $(A,B)$ is stabilizable, and $(A,C)$ is detectable.
\end{assumption}

\subsection{Optimal control problem}
We will employ a finite-horizon control policy with input at time $k$ in the form
\[
u_{i|k} = \kappa_i(\theta_k ,\U_{k+i-1},\I_{k+i})
\]
where $u_{i|k}$ for $i=0,1,\ldots$ is the prediction of $u_{k+i}$ at time $k$, and $\theta_k$ is a vector of controller parameters at time $k$. The dependence of $\kappa_i(\cdot)$ on the sets $\U_{k+i-1}$ and $\I_{k+i}$ ensures causality and the dependence on $\theta_k$ is chosen so that the optimal parameter vector, denoted $\theta_k^\ast$, will be the solution of a convex problem.
%ensure convexity of an optimization problem that is solved at time $k$ to determine the optimal parameter $\theta_k^\ast$.

\begin{assumption}\label{assumption:info}
(i). The probability, $\lambda$, of successfully receiving sensor measurements is known. (ii). When $\theta_k^\ast$ is computed, $(z_{k+i},\gamma_{k+i})$ are unknown for all $i\geq 0$.
\end{assumption}

Assumption~\ref{assumption:info} requires $\theta_k^\ast$ to be a function of $\U_{k-1}$ and $\I_{k-1}$, and we therefore assume that $\theta_k^\ast$ is computed online prior to the $k$th sampling instant. However $(z_k, \gamma_k)$ is known when the control law
\[
u_k = \kappa_0(\theta_k^\ast , \U_{k-1}, \I_{k})
\]
is applied to the plant.
%the arrival of the sensor data $z_{k}$, $\gamma_{k}$ at the controller.
%

We consider the problem of minimizing the discounted sum of expected future values of ${\|x_{k}\|_Q^2 \!+\! \|u_{k}\|_R^2}$, where $Q\!\succeq\! 0$ and $R\!\succ\! 0$. This minimization is subject to a constraint on the discounted sum of second moments of an auxiliary output, defined for given matrix $H$  by $\xi_k \!=\! H x_{k}$, so that
%(with $\PP_k\{\cdot\}$ denoting $\PP\{\cdot \,|\, \U_{k-1},\I_{k-1}\}$)
%\begin{equation}\label{eq:orig_opt_problem}
\vspace{-1.5mm}
\begin{align}
\theta_k^\ast = \arg\min_{\theta_k} \
& \sum_{i=0}^{\infty} \beta^i \EE_k\{  \|x_{i|k}\|_Q^2 + \|u_{i|k}\|_R^2 \}
\nonumber \\
\text{s.t.} \
& \sum_{i=0}^\infty \beta^i \EE_k\{ \| H x_{i|k}\|^2 \} \leq \epsilon .
\label{eq:orig_opt_problem}
\end{align}
%\end{equation}
%the probabilities and expectations are conditioned on $\U_{k-1}$, $\Z_{k-1}$, $\Gamma_{k-1}$ and
Here $\beta\in(0,1)$ is a discounting factor and $\epsilon$ is a given bound on this infinite discounted sum of second moments.
%and the value of $\mu_{k}$ is assumed to be chosen in order to guarantee recursive feasibility of (\ref{eq:opt}) while ensuring that the closed loop system satisfies the original constraints of the problem.
%
%We use sufficient conditions based on Chebyshev's inequality to bound the cost of~(\ref{eq:orig_opt_problem}).
%This will allow~(\ref{eq:orig_opt_problem}) to be formulated as a convex problem, as required for computational tractability, while providing robustness to assumptions on the distribution of disturbances and sensor noise.
%Using $\PP\{ \|x_k\|_Q^2 + \|u_k\|_R^2 \leq 1\} \leq \EE\{ \|x_k\|_Q^2 + \|u_k\|_R^2\}$ \citep[Sec.~V.7]{feller71}, 
Instead of solving \eqref{eq:orig_opt_problem} directly, the control problem to be solved at time $k$ is given by
\vspace{-1mm}
\begin{equation}\label{eq:opt}
\begin{aligned}
\theta_k^\ast = \arg\min_{\theta_k} \
& \sum_{i=0}^{\infty} \beta^i \EE_k\{  \|x_{i|k}\|_Q^2 + \|u_{i|k}\|_R^2 \}  \\
\text{s.t.} \
& \sum_{i=0}^\infty \beta^i \EE_k\{ \| H x_{i|k} \|^2\} \leq \mu_{k} .
\end{aligned}
\end{equation}
Here $\mu_0=\epsilon$ and, for all $k>0$, $\mu_k$ is chosen as described in Section~\ref{sec:properties} to ensure that (\ref{eq:opt}) is recursively feasible and that the constraint in (\ref{eq:orig_opt_problem}) is satisfied with $k=0$ by the closed loop system.

\section{Controller parameterization}
\label{sec:control_parameterization}

Consider the output feedback control law defined by an observer and an affine feedback law:
\begin{subequations}\label{eq:basic_control_law}
\begin{gather}
\xhat_k = A \tilde{x}_{k-1} + Bu_{k-1},~
\tilde{x}_k = \xhat_{k} + \gamma_k M (y_k - C\xhat_k), \\
u_k = K \tilde{x}_k  + c_k.
\end{gather}
\end{subequations}
with $\xhat_0 \!=\! \EE\{x_0\}$, where $\xhat_k$ and $\tilde{x}_k$ are the a priori estimate and the posteriori estimate of $x_k$, respectively. 
A simplistic parameterization of the predicted control law $\kappa_i(\cdot)$ could be obtained if the observer gain $M$ and feedback gain $K$ were fixed and the optimization variables in problem (\ref{eq:opt}) were defined as $\theta_k \!=\! \{c_{0|k},\ldots,c_{N-1|k}\}$ for some fixed $N$, with the predicted control sequence defined as $u_{i|k} \!=\! K \tilde{x}_{i|k}  \!+\! c_{i|k}$.
Although this would  require a number of optimization variables that grows only linearly with $N$, the parameters $\{c_{0|k},\ldots,c_{N-1|k}\}$ constitute an open loop control sequence that does not vary with the future measurement noise and disturbance realizations. This is likely to provide poor performance and small sets of feasible initial conditions when the probability of packet loss is non-zero.

By using a parameterization that allows the dependence of the predicted control sequence on future realizations of model uncertainty to be optimized, the predicted probability distributions of states and control inputs can be controlled explicitly. This provides flexibility to balance conflicting requirements for performance and constraint satisfaction.
However, similarly to the case of predicted control laws in which state feedback gains are decision variables~\citep{lofberg03,goulart06}, the cost and constraints of problem~(\ref{eq:opt}) are nonconvex if time-varying gains $M$, $K$ are considered as optimization variables.
On the other hand, if predicted control inputs are parameterized in terms of affine functions of the future output measurements received by the controller, then the dependence of the first and second moments of predicted states and inputs on controller parameters is convex.
Moreover, by incorporating affine terms in the future innovation sequence, a predicted control law with arbitrary linear dependence of $\kappa_i(\cdot)$ on the received sensor measurements can be obtained.
This approach allows the future control sequence to be optimized at every sampling instant, including those at which information from sensors is lost.

We therefore express the $i$ steps ahead predicted control input $u_{i|k}$, for all $i = 0,1,\ldots$, as
\begingroup
\setlength{\abovedisplayskip}{-1pt}
\setlength{\belowdisplayskip}{0pt}
\begin{subequations}\label{eq:predicted_control_law}
\begin{align}
u_{i|k} &= K \xhat_{i|k} + c_{i|k} + d_{i|k},\\
d_{i|k} &= \gamma_{0|k}L_{i,0|k} (y_{0|k} \!-\! C\xhat_{0|k})  + \gamma_{1|k}L_{i,1|k} (y_{1|k} \!-\! C\xhat_{1|k})  \nonumber \\
&\quad + \cdots +
\gamma_{i|k}L_{i,i|k} (y_{i|k} - C\xhat_{i|k}),
\\
\xhat_{i+1|k} &= A\xhat_{i|k} + B u_{i|k} + \gamma_{i|k} AM (y_{i|k} - C\xhat_{i|k} ) \label{eq:predicted_xhat}
\end{align}
\end{subequations}
\endgroup
where $c_{i|k} = 0$ and $L_{i,j|k} = 0$ for all $i\geq N$. Here $\gamma_{i|k}$ and $y_{i|k}$ are random variables, denoting the $i$-step-ahead predicted packet loss and sensor measurement at time $k$, respectively.
Then, for all $i = 0,1,\ldots$ the predicted state estimate satisfies
\begin{equation}\label{eq:nominal_state_estimate}
\xhat_{i+1|k} = \Phi \xhat_{i|k} + B(c_{i|k} + d_{i|k}) + \gamma_{i|k} A M (y_{i|k} - C \xhat_{i|k})
\end{equation}
where $\Phi \defeq A+BK$. Since $x_{i+1|k} = A x_{i|k} + Bu_{i|k} + Dw_{i|k}$ the predicted estimation error evolves according to
\begin{equation}\label{eq:state_prediction}
x_{i+1|k} - \xhat_{i+1|k} = \Psi_{i|k}(x_{i|k} - \xhat_{i|k})
- \gamma_{i|k}AMv_{i|k} + D w_{i|k}
\end{equation}
with $\Psi_{i|k} \defeq A (I - \gamma_{i|k}M C)$.
These relationships allow the first and second moments of $x_{i|k}$ to be determined in terms of the decision variable $\theta_k$, which consists of the parameters
$\{c_{0|k},\ldots,c_{N-1|k}\}$
and feedback gains
$L_{0,0|k}$, $\{L_{1,0|k},L_{1,1|k}\}$, $\ldots, \{L_{N-1,0|k},\ldots,L_{N-1,N-1|k}\}$.

The gains $K$ and $M$ in the predicted control law (\ref{eq:predicted_control_law}a-c) are chosen of\mbox{}f\mbox{}line and satisfy the following assumption.

\begin{assumption}\label{assumption:stabilizing_gains}
$\xi_{i+1} = (A+BK) \xi_i$ is asymptotically stable and $\xi_{i+1} = A(I - \gamma_iMC) \xi_i$ is mean-square stable \citep{Kushner1971}.
\end{assumption}

\begin{rem}
Gains $K$ and $M$ exist satisfying Assumption~\ref{assumption:stabilizing_gains} if Assumption~\ref{assumption:ctrb_and_obsv} holds and if the  probability, $\lambda$, of successfully receiving a sensor measurement is greater than some critical value~\citep[e.g.][]{sinopoli04}. Suitable choices for $K$, $M$
are the optimal gains for (\ref{eq:opt}) in the absence of constraints, or the certainty equivalent LQ feedback gain for a problem with state and control weighting matrices $Q$ and $R$ and the steady state Kalman filter gain~\citep{sinopoli04}. We note also that time-varying gains $K_k$, $M_k$ can be used within the framework of this paper, provided their dependence on $\gamma_k$ is known in advance.
\end{rem}

\subsection{Controller implementation}
\label{sec:control_implementation}

The control law is implemented by the following procedure.
\begin{enumerate}[(i).]
\item
Given $\mathcal{U}_{k-1}$ and $\I_{k-1}$, solve problem (\ref{eq:opt}) for $\theta^\ast_k$.
\item
Given $\gamma_k$ and $z_k=\gamma_k y_k$:
\begin{enumerate}[(a).]
\item
apply the control input
\begingroup
\setlength{\abovedisplayskip}{-0.5pt}
\setlength{\belowdisplayskip}{-0.5pt}
\[
u_k = K\xhat_k + c_{0|k}^\ast + \gamma_k L^\ast_{0,0|k}(y_k - C\xhat_k),
\]
\endgroup
\item
update the state estimate
\begingroup
\setlength{\abovedisplayskip}{-0.5pt}
\setlength{\belowdisplayskip}{-1.5pt}
\[
\xhat_{k+1} = A\xhat_k + B u_k + \gamma_k A M(y_k - C\xhat_k).
\]
\endgroup
\end{enumerate}
\end{enumerate}

Note that this receding horizon control law includes (\ref{eq:basic_control_law}) as a special case, since $u_{k}$ and $\xhat_{k+1}$ in step (ii) would be equal to their counterparts in (\ref{eq:basic_control_law}) if $(c_{0|k}^\ast,L_{0,0|k}^\ast) = (c_k,KM)$.%
%, and where $M$ denotes the steady state optimal estimator gain when constraints are inactive.
%

\section{Predicted state and control sequences}
\label{sec:MPC_predictions}
%Given a sequence $\{ x_{0|k},\ldots, x_{N|k} \}$, let $\bx_k$ denote the vector $[x_{0|k}^\top \ \cdots \ x_{N|k}^\top]^\top$, and,
To simplify notation we express the predicted control law in terms of vectorized sequences, with
$\bx_k$ denoting the vectorized state sequence $\{x_{i|k}\}_{i=0}^{N-1}$,
$\bxhat_k$ the estimate sequence $\{\xhat_{i|k}\}_{i=0}^{N-1}$,
 $\bu_k$ the predicted control sequence $\{u_{i|k}\}_{i=0}^{N-1}$,
${\bf c}_k$ the predicted control perturbations $\{c_{i|k}\}_{i=0}^{N-1}$,
$\bw_k$ the disturbances $\{w_{i|k}\}_{i=0}^{N-1}$,
$\bv_k$ the sensor noise sequence $\{v_{i|k}\}_{i=0}^{N-1}$, and
$\bs{\zeta}_k$ the future innovation sequence $\{\gamma_{i|k} (y_{i|k}-C\xhat_{i|k}) \}_{i=0}^{N-1}$ at time $k$.
For a given sequence of matrices $\{\Psi_{i|k}\}_{i=0}^{N-1}$ and matrix $B$ let
\vspace{-1mm}
\begin{alignat*}{2}
{\bf S}_{\Psi} &= \begin{bmatrix} I \\ \Psi_{0|k} \\ \vdots \\ \prod_{i=N-2}^0 \Psi_{i|k} \end{bmatrix}\!,
&~
{\bf T}_{(\Psi,B)} &=
\begin{bmatrix}
0 & \cdots & & 0 \\
B & & &\\
\vdots & &\ddots & \\
\prod_{i=N-2}^1 \Psi_{i|k} B & \cdots & B & 0
\end{bmatrix},
\\
S^N_{\Psi} &= \prod_{i=N-1}^0 \Psi_{i|k},
&~
T^N_{(\Psi,B)} &= \begin{bmatrix} \prod_{i=N-1}^1 \Psi_{i|k} B & \cdots & B \end{bmatrix}
% {\bf S}_{\Psi} = \begin{bmatrix} I \\ A \\ \vdots \\ A^{N} \end{bmatrix},
% \quad
% {\bf T}_{(A,B)} = \begin{bmatrix} 0 & \cdots & 0 \\ B & & \\ \vdots & \ddots & \\ A^{N-1}B & \cdots & B \end{bmatrix}
\end{alignat*}
where $\prod_{i=m}^n \Psi_{i|k} = \Psi_{m|k} \cdots \Psi_{n|k}$ for $m \geq n$,
and define
\begin{align*}
{\bf L}_k &= \begin{bmatrix} L_{0,0|k} & & & \\ L_{1,0|k} & L_{1,1|k} & & \\ \vdots & \vdots & \ddots & \\ L_{N-1,0|k} & L_{N-1,1|k} & \cdots & L_{N-1,N-1|k} \end{bmatrix} ,\\
\bs{\Gamma}_k &= \diag\{ \gamma_{0|k}, \ldots, \gamma_{N-1|k} \} \otimes I_{n_y \times n_y} ,
\end{align*}
${\bf K} = I_{N\times N}\otimes K$, ${\bf M} = I_{N\times N}\otimes M$ and ${\bf C} = I_{N\times N} \otimes C$.
Then from (\ref{eq:state_prediction}) we have
\begin{equation}\label{eq:pred_xdiff}
\bx_k - \bxhat_k = {\bf S}_{\Psi} (x_k - \xhat_k)
- {\bf T}_{(\Psi,A)} {\bf M}\bs{\Gamma}_k \bv_k
+ {\bf T}_{(\Psi,D)} \bw_k
% \\
% x_{N|k} - \xhat_{N|k} &= S^N_{\Psi} (x_k - \xhat_k)
% - T^N_{(\Psi,A)} {\bf M}\bs{\Gamma}_k \bv_k
% + T^N_{(\Psi,D)} \bw_k ,
\end{equation}
while (\ref{eq:nominal_state_estimate}) and (\ref{eq:predicted_control_law}b) give
\[
\bxhat_k = {\bf S}_{\Phi} \xhat_k + {\bf T}_{(\Phi,B)}({\bf c}_k + {\bf L}_k\bs{\zeta}_k) + {\bf T}_{(\Phi,A)} {\bf M} \bs{\zeta}_k .
\]
Here $\bs{\zeta}_k = \bs{\Gamma}_k {\bf C}(\bx_k - \bxhat_k) + \bs{\Gamma}_k\bv_k$ and ${\bf S}_{\Phi}$, ${\bf T}_{(\Phi,B)}$ are defined (analogously to ${\bf S}_{\Psi}$, ${\bf T}_{(\Psi,B)}$) in terms of $\Phi$ and $B$. Hence
\begingroup
\setlength{\abovedisplayskip}{-10pt}
\setlength{\belowdisplayskip}{-0.5pt}
\begin{subequations}\label{eq:pred_xhat_u}
\begin{align}
\bxhat_k &= {\bf S}_{\Phi}\xhat_k \!+ {\bf T}_{(\Phi,B)}{\bf c}_k \!+ ({\bf T}_{(\Phi,B)}{\bf L}_k + {\bf T}_{(\Phi,A)}{\bf M}) \bs{\zeta}_k,
\\
\bu_k &= {\bf K} \bxhat_k + {\bf c}_k + {\bf L}_k \bs{\zeta}_k
.
% \xhat_{N|k} &= S^N_{\Phi}\xhat_k
% + T^N_{(\Phi,B)}{\bf c}_k + (T^N_{(\Phi,B)}{\bf L}_k + T^N_{(\Phi,A)}{\bf M}) \bs{\Gamma}_k \bs{\zeta}_k \\
\end{align}
\end{subequations}
\endgroup
Clearly the predicted estimation error, state and control sequences in (\ref{eq:pred_xdiff}) and (\ref{eq:pred_xhat_u}a,b) depend linearly on the decision variables
$\theta_k \defeq ({\bf c}_k, {\bf L}_k)$.

\subsection{First and second moments of predicted sequences}

In order to express the cost and constraints of problem~(\ref{eq:opt}) in terms of the parameterization introduced in Section~\ref{sec:control_parameterization}, we derive in this section expressions for the means and variances of predicted state and control sequences.
%
% In the following all expectations are conditional on $\U_{k-1},\Z_{k-1},\Gamma_{k-1}$, and the realisations of $\gamma_k$ and $z_k=\gamma_ky_k$ are assumed unknown.
% %S

First consider the state $x_k$ of the plant (\ref{eq:system}a) and the state estimate update $\xhat_k$ in step (ii) of the controller implementation in Section~\ref{sec:control_implementation}.
By assumption we have $\EE\{x_0\} = \xhat_0$ and $\EE\{w_{k}\}=0$, $\EE\{v_{k}\}=0$ for all $k\geq 0$, and hence the update of state estimates $\xhat_k$ in step (ii)(b) ensures that
\vspace{-1.5mm}
\begin{equation}\label{eq:cl_meanx}
\EE_k\{x_k \} = \xhat_k
\end{equation}
for all $k\geq 1$. Furthermore, from (\ref{eq:system}a) we have
\begingroup
\setlength{\abovedisplayskip}{-0.8pt}
\setlength{\belowdisplayskip}{-0.5pt}
\[
x_k - \xhat_k = \Psi_{k-1} (x_{k-1} - \xhat_{k-1}) - \gamma_{k-1} AM v_{k-1} + D w_{k-1}
\]
\endgroup
for all $k\geq 1$. Let $\Sigma_k$ denote the second moment of the state estimate error at time $k$:
\[
\Sigma_k \defeq \EE_k \bigl\{(x_k-\xhat_k)(x_k - \xhat_k)^\top \bigr\}.
\]
Then $\Sigma_k$ evolves according to
\begin{equation}\label{eq:cl_varx}
\Sigma_k \!=\! \Psi_{k-1} \Sigma_{k-1} \Psi_{k-1}^\top \!+ \gamma_{k-1}AM\Sigma_v M^\top A^\top \!+ D \Sigma_w D^\top
\end{equation}
for all $k\geq 1$, with initial condition $\Sigma_0$, and by Assumption \ref{assumption:stabilizing_gains} $\EE \{\Sigma_k\}$ remains upper bounded $\forall k$.

We first derive the first and second moments of the predicted state sequence $\bx_k$ and control sequence $\bu_k$:

\begin{prop}\label{prop:moments_x_u}
Let $\bs{\pi}_k$, $\bs{\Pi}_k$, and $\bs{\Omega}_k$ be defined
\begingroup
\setlength{\abovedisplayskip}{-0.5pt}
\setlength{\belowdisplayskip}{-0.5pt}
\begin{gather*}
\bs{\pi}_k = {\bf S}_{\Phi} \xhat_k + {\bf T}_{(\Phi,B)}{\bf c}_k ,
\quad
\bs{\Pi}_k = {\bf T}_{(\Phi,B)} {\bf L}_k + {\bf T}_{(\Phi,A)}{\bf M} ,
\\
\bs{\Omega}_k = \EE_k \biggl\{
\begin{bmatrix} \bx_k - \bxhat_k \\ \bs{\zeta}_k \end{bmatrix}
\begin{bmatrix} \bx_k - \bxhat_k \\ \bs{\zeta}_k \end{bmatrix}^\top
\biggr\}.
\end{gather*}
\endgroup
Then
\vspace{-2mm}
\begin{subequations}\label{eq:1st_pred_x_u}
\begin{align}
&\EE_k \{\bx_k\}  = \EE_k\{\bxhat_k\} = \bs{\pi}_k,
\\
&\EE_k \{\bu_k\} = {\bf K} \EE_k\{\bxhat_k\} + {\bf c}_k
= {\bf K}\bs{\pi}_k + {\bf c}_k,
\end{align}
\end{subequations}
and
\vspace{-2mm}
\begin{subequations}\label{eq:2nd_pred_x_u}
\begin{align}
&\EE_k\{\bx_k\bx_k^\top\} = \bs{\pi}_k\bs{\pi}_k^\top + \begin{bmatrix} I & \bs{\Pi}_k\end{bmatrix} \bs{\Omega}_k \begin{bmatrix} I \\ \bs{\Pi}_k^\top\end{bmatrix},
\\
&\EE_k\{\bu_k\bu_k^\top\} = ({\bf K}\bs{\pi}_k + {\bf c}_k)({\bf K}\bs{\pi}_k + {\bf c}_k)^\top
\nonumber\\
&\qquad\qquad\quad + \begin{bmatrix} 0 & {\bf L}_k+{\bf K}\bs{\Pi}_k\end{bmatrix} \bs{\Omega}_k \begin{bmatrix} 0 \\ ({\bf L}_k+{\bf K}\bs{\Pi}_k)^\top\end{bmatrix} .
\end{align}
\end{subequations}

%where all expectations are conditioned on $\U_{k-1}$, $\I_{k-1}$.
\end{prop}
\vspace{-1mm}
\bpf
From~(\ref{eq:pred_xdiff}), (\ref{eq:cl_meanx}) we have $\EE_k \{\bx_k - \bxhat_k \} = 0$. Therefore $\bs{\zeta}_k = \bs{\Gamma}_k{\bf C}(\bx_k - \bxhat_k) + \bs{\Gamma}_k{\bf v}_k$ implies $\EE_k\{\bs{\zeta}_k \} = 0$
and (\ref{eq:1st_pred_x_u}a,b) follow from the expectations of (\ref{eq:pred_xhat_u}a,b).
To determine the second moments of $\bx_k$ and $\bu_k$, let
\vspace{-1mm}
\[
{\bf X}_k \defeq \EE_k\biggl\{
\begin{bmatrix} \bx_k - \bxhat_k \\ \bxhat_k \end{bmatrix}
\begin{bmatrix} \bx_k - \bxhat_k \\ \bxhat_k \end{bmatrix}^\top
\biggr\}\ .
\]
Then from (\ref{eq:pred_xdiff}) and (\ref{eq:pred_xhat_u}a) we have
\vspace{-1mm}
\begin{equation}\label{eq:x_covar}
{\bf X}_k =
\begin{bmatrix} 0 & 0 \\ 0 & \bs{\pi}_k \bs{\pi}_k^\top \end{bmatrix}
+
\begin{bmatrix} I & 0 \\ 0 & \bs{\Pi}_k\end{bmatrix} \bs{\Omega}_k \begin{bmatrix} I & 0 \\ 0 & \bs{\Pi}_k\end{bmatrix}^\top ,
\end{equation}
and (\ref{eq:2nd_pred_x_u}a,b) follow from
$\EE_k\{\bx_k \bx_k^\top \} = \begin{bmatrix} I & I \end{bmatrix} {\bf X}_k \begin{bmatrix} I & I \end{bmatrix}^\top$
and (\ref{eq:pred_xhat_u}a,b), respectively.
\epf

Since $\bs{\pi}_k$ and $\bs{\Pi}_k$ are linear in $({\bf c}_k,{\bf L}_k)$ and $\bs{\Omega}_k$ is independent of $({\bf c}_k,{\bf L}_k)$,
it is clear from (\ref{eq:1st_pred_x_u}a,b) and (\ref{eq:2nd_pred_x_u}a,b)
that the first moments of the predicted state and input sequences are linear in $\theta_k = ({\bf c}_k,{\bf L}_k)$ while their second moments are quadratic functions of $\theta_k$.

To determine $\bs{\Omega}_k$, note that $\bx_k - \bxhat_k$ and $\bs{\zeta}_k$ can be written
\vspace{-3mm}
\[
\bx_k - \bxhat_k = F(\bs{\Gamma}_k)q_k ,
\quad
\bs{\zeta}_k = G(\bs{\Gamma}_k)q_k  ,
\quad
q_k = \begin{bmatrix} x_k - \hat{x}_k \\ {\bf v}_k \\ {\bf w}_k \end{bmatrix},
\]
with
$F(\bs{\Gamma}_k) = \begin{bmatrix} {\bf S}_{\Psi} & {-{\bf T}_{(\Psi,A)}}{\bf M} \bs{\Gamma}_k & {\bf T}_{(\Psi,D)} \end{bmatrix}$ and
$G(\bs{\Gamma}_k) = \bs{\Gamma}_k{\bf C}F(\bs{\Gamma}_k) + \begin{bmatrix} 0 & \bs{\Gamma}_k & 0 \end{bmatrix}$. So, by the law of total expectation,
\vspace{-2mm}
\begin{equation}\label{eq:cov_mat}
\bs{\Omega}_k \!=\!\!
\sum_j \!\begin{bmatrix} F(\bs{\Gamma}^{(j)}) \\ G(\bs{\Gamma}^{(j)}) \end{bmatrix} \!\EE\{ q_k q_k^\top \}\! \begin{bmatrix} F(\bs{\Gamma}^{(j)}) \\ G(\bs{\Gamma}^{(j)}) \end{bmatrix}^\top \!\!\PP \{ \bs{\Gamma}_k \!=\! \bs{\Gamma}^{(j)} \} ,
\end{equation}
where $\EE\{ q_k q_k^\top \}$ is the block-diagonal matrix:
\begin{gather*}
\EE\{ q_k q_k^\top \} = \diag\{ \Sigma_k, \bar{\Sigma}_v, \bar{\Sigma}_w\},
\\
\bar{\Sigma}_v =I_{N \times N} \otimes \Sigma_v,
\quad
\bar{\Sigma}_w =I_{N \times N} \otimes \Sigma_w,
\end{gather*}
and where $\bs{\Gamma}^{(j)}$ for $j =1,\ldots,2^N$ enumerates the $2^N$ matrices with binary-valued diagonal elements defined by
\begin{align*}
&\bs{\Gamma}^{(1)} = 0, \quad
\bs{\Gamma}^{(2)} = \diag\{0,\ldots,0,1\}\otimes I_{n_y\times n_y}  \quad \ldots
\\
&
\ldots \quad \bs{\Gamma}^{(2^N-1)} = \diag\{1,\ldots,1 ,0\}\otimes I_{n_y\times n_y}, \quad
\bs{\Gamma}^{(2^N)} = I .
\end{align*}

\begin{rem}
$\bs{\Omega}_k $ in (\ref{eq:cov_mat}) can be computed conveniently via
\begin{align*}
&\vvec(\bs{\Omega}_k) =
\\
&\biggl( \sum_j\!
\begin{bmatrix} F(\bs{\Gamma}^{(j)}) \\ G(\bs{\Gamma}^{(j)}) \end{bmatrix} \!\otimes\!
\begin{bmatrix} F(\bs{\Gamma}^{(j)}) \\ G(\bs{\Gamma}^{(j)}) \end{bmatrix}
\PP \{ \bs{\Gamma}_k = \bs{\Gamma}^{(j)} \} \!\biggr)
\!\vvec\Bigl(\Bigl[\begin{smallmatrix} \!\Sigma_k\! & & \\ & \!\bar{\Sigma}_v\! & \\ & & \!\bar{\Sigma}_w\! \end{smallmatrix} \Bigr]\Bigr)
\end{align*}
where the first term on the RHS can be determined offline given the probability distribution of $\gamma_k$. This allows $\bs{\Omega}_k$ to be computed online using the current value of $\Sigma_k$ with a single matrix-vector multiplication.
\end{rem}

% Similarly ${\bf U}_k = \EE\{\bu_k\bu_k^\top\}$ satisfies
% \begin{equation}\label{eq:u_covar}
% {\bf U}_k =
% \bigl( {\bf L}_k + \bs{\Pi})
% %\bigl( (I + {\bf T}_{(\Phi,B)}) {\bf L}_k + {\bf T}_{(\Phi,A)}{\bf M}\bigr)
% \begin{bmatrix} 0 & 0 \\ 0 & \bs{\Omega}_k \end{bmatrix}
% \bigl( {\bf L}_k + \bs{\Pi})^\top
% %\bigl( (I + {\bf T}_{(\Phi,B)}) {\bf L}_k + {\bf T}_{(\Phi,A)}{\bf M}\bigr)^\top
% + ( {\bf K} \bs{\pi}_k + {\bf c}_k) ( {\bf K} \bs{\pi}_k + {\bf c}_k)^\top
%  % \bigl( {\bf K} {\bf S}_{\Phi} \xhat_k + (I + {\bf K}{\bf T}_{(\Phi,B)}){\bf c}_k\bigr) \bigl( {\bf K} {\bf S}_{\Phi} \xhat_k + (I + {\bf K}{\bf T}_{(\Phi,B)}){\bf c}_k\bigr)^\top ,
% \end{equation}
% and, defining $E_i = [e_i^\top \ \ e_i^\top] E$ where $e_i$ is the $i$th column of $I_{N\times N}$, we have
% \begin{equation}\label{eq:output_covar}
% \EE\{(E x_{i|k} )^\top (E x_{i|k}) \} = \tr (E_i^\top E_i {\bf X}_k)  .
% \end{equation}

Using the same arguments as the proof of Proposition~\ref{prop:moments_x_u},
it can be verified that
\vspace{-1mm}
\begin{align}
X_{N|k} &=
\EE_k \biggl\{\begin{bmatrix} x_{N|k} - \xhat_{N|k} \\ \xhat_{N|k} \end{bmatrix}
\begin{bmatrix} x_{N|k} - \xhat_{N|k} \\ \xhat_{N|k} \end{bmatrix}^\top
\biggr\} \nonumber \\
&=
\begin{bmatrix} 0 & 0 \\ 0 & \pi_{N|k} \pi_{N|k}^\top \end{bmatrix}
+
\begin{bmatrix} I & 0 \\ 0 & \Pi_{N|k}\end{bmatrix} \Omega_{N|k} \begin{bmatrix} I & 0 \\ 0 & \Pi_{N|k}\end{bmatrix}^\top
\label{eq:xN_covar}
\end{align}
where
\vspace{-1mm}
\begin{align*}
\Pi_{N|k} &= T^N_{(\Phi,B)} {\bf L}_k + T^N_{(\Phi,A)} {\bf M},~ 
\pi_{N|k} = S^N_{\Phi} \xhat_k + T^N_{(\Phi,B)} {\bf c}_k, \\
\Omega_{N|k} &= \sum_j\! \begin{bmatrix} F_N(\bs{\Gamma}^{(j)}) \\ G (\bs{\Gamma}^{(j)}) \end{bmatrix}
\!\Bigl[ \begin{smallmatrix}
\!\Sigma_k  & & \\ & \!\bar{\Sigma}_v\! & \\  & & \bar{\Sigma}_w\!
\end{smallmatrix} \Bigr]\!
\begin{bmatrix} F_N (\bs{\Gamma}^{(j)}) \\ G (\bs{\Gamma}^{(j)}) \end{bmatrix}^{\!\top}\!\! \PP \{ \bs{\Gamma}_k \!=\! \bs{\Gamma}^{(j)} \}
\end{align*}
with
$F_N(\bs{\Gamma}_k) = \begin{bmatrix} S^N_{\Psi} & -T^N_{(\Psi,A)}{\bf M} \bs{\Gamma}_k & T^N_{(\Psi,D)} \end{bmatrix}$.
%  \\
% and $G_N(\bs{\Gamma}_k) = \gamma_{N|k} C F_N(\bs{\Gamma}_k) + \gamma_{N|k}\begin{bmatrix} 0 & e_N^\top\otimes I_{n_y\times n_y} & 0 \end{bmatrix} $.

\section{Cost and constraints} \label{sec:cost and constraints}

We next show that the cost and constraints of (\ref{eq:opt}) can be expressed as convex functions of $\theta_k=({\bf c}_k,{\bf L}_k)$.
First note that the objective in (\ref{eq:opt}) can be written
\vspace{-1mm}
\begin{multline}
\sum_{i=0}^{\infty} \beta^i \EE_k \bigl\{ \| x_{i|k} \|_Q^2 + \| u_{i|k}\|_R^2 \bigr\}
\\
=
\tr (\Q_\beta{\bf X}_k) + \tr(\R_\beta{\bf U}_k) + f_N(\theta_k,\xhat_k,\Sigma_k)
\label{eq:cost_bound}
\end{multline}
where ${\bf X}_k$ is given by (\ref{eq:x_covar}), and
\begin{align*}
&{\bf U}_k \defeq \EE_k \{\bu_k\bu_k^\top \},
\Q_\beta \! \defeq\!
\ones_{2\times 2} \otimes\diag\{ Q , \beta Q, \ldots, \beta^{N-1}Q \} ,
\\
&\R_\beta  \defeq \diag\{ R , \beta R, \ldots, \beta^{N-1}R \} ,
\\
&f_N(\theta_k, \xhat_k,\Sigma_k) \defeq \sum_{i=N}^{\infty} \beta^i \EE_k \bigl\{ \| x_{i|k} \|_Q^2 + \| u_{i|k}\|_R^2 \bigr\} .
\end{align*}
Since $\Q_\beta\succeq 0$ and $\R_\beta \succ 0$, the term $\tr (\Q_\beta{\bf X}_k) + \tr(\R_\beta{\bf U}_k)$ in (\ref{eq:cost_bound}) can be expressed as a convex quadratic function of $\theta_k = ({\bf c}_k,{\bf L}_k)$ using (\ref{eq:2nd_pred_x_u}b) and (\ref{eq:x_covar}).
%
% Then the cost (\ref{eq:cost_bound}) can be expressed as a convex quadratic function of $\theta_k = ({\bf c}_k,{\bf L}_k)$ by using (\ref{eq:2nd_pred_x_u}b) and (\ref{eq:x_covar}) to evaluate $\tr (\Q{\bf X}_k) + \tr(\R{\bf U}_k)$.
%
To determine the terminal term, $f_N(\theta,\xhat_k,\Sigma_k)$, let $P_k = \sum_{i=N}^\infty \beta^i X_{i|k}$, where
\vspace{-2mm}
\[
X_{i|k} = \EE_k \biggl\{\begin{bmatrix} x_{i|k} - \xhat_{i|k} \\ \xhat_{i|k} \end{bmatrix}
\begin{bmatrix} x_{i|k} - \xhat_{i|k} \\ \xhat_{i|k} \end{bmatrix}^\top\biggr\} .
\]
Then for $i\geq N$ we have
\[
X_{i+1|k} = \EE \bigl\{ \tilde{\Psi}(\gamma) X_{i|k} \tilde{\Psi}^\top(\gamma) \bigr\}
+ \EE \bigl\{ \tilde{D}(\gamma) \bigl[\begin{smallmatrix}\Sigma_v & \\ & \Sigma_w\end{smallmatrix}\bigr]
\tilde{D}^\top (\gamma) \bigr\}
\]
where
\vspace{-2mm}
\[
\tilde{\Psi}(\gamma) = \begin{bmatrix} A(I - \gamma MC) & 0 \\ \gamma A M C & \Phi \end{bmatrix} ,
\quad
\tilde{D}(\gamma) = \begin{bmatrix} -\gamma A M & D \\ \gamma A M & 0\end{bmatrix} ,
\]
and $\gamma$ is a random variable identically distributed as $\gamma_k$.
Hence
%$\EE\bigl\{ \tilde{\Psi}(\gamma) P_k \tilde{\Psi}^\top(\gamma) \bigr\} = \sum_{i=N}^\infty \beta^i \EE\bigl\{ \tilde{\Psi}(\gamma) X_{i|k} \tilde{\Psi}^\top(\gamma) \bigr\}$ and
\begin{align*}
& \EE\bigl\{ \tilde{\Psi}(\gamma) P_k \tilde{\Psi}^\top(\gamma) \bigr\} \\
&=
\sum_{i=N}^\infty \beta^i \bigl( X_{i+1|k} - \EE\{\tilde{D}(\gamma) \bigl[\begin{smallmatrix}\Sigma_v & \\ & \Sigma_w\end{smallmatrix}\bigr]
\tilde{D}^\top (\gamma) \}\bigr)
\\
&= \beta^{-1} ( P_k \!-\! \beta^N X_{N|k}) -\!\frac{\beta^N}{1-\beta} \EE \{\tilde{D}(\gamma)
\bigl[\begin{smallmatrix}\Sigma_v & \\ & \Sigma_w\end{smallmatrix}\bigr]
\tilde{D}^\top (\gamma) \} ,
\end{align*}
and the terminal term $f_N(\theta_k,\xhat_k,\Sigma_k)$ in (\ref{eq:cost_bound}) is equal to
\[
%f_N(\theta_k,\xhat_k,\Sigma_k) =
\tr\Bigl( \Bigl[\begin{smallmatrix}  Q & ~Q \\ Q & ~Q + K^\top\! RK \end{smallmatrix}\Bigr] P_k \Bigr)
\]
with the additional constraint
% It follows that the problem of minimizing
% \[
% \sum_{i=0}^{\infty} \beta^i \bigl( \EE \{ \| x_{i|k} \|_Q^2 \} + \EE\{ \| u_{i|k}\|_R^2 \} \bigr)
% \]
% over $\theta_k$ is equivalent to the problem of minimizing
% \begin{equation}\label{eq:infinite_horizon_cost}
% \tr (\Q{\bf X}_k) + \tr(\R{\bf U}_k) +
% \tr\Bigl( \begin{bmatrix} Q & Q \\ Q & Q + K^\top R K\end{bmatrix} P_k\Bigr)
% \end{equation}
% over $\theta_k$, $P_k$ subject to (\ref{eq:x_covar}), (\ref{eq:u_covar}), (\ref{eq:xN_covar}) and
\begin{align}
P_k &\succeq \beta \EE\bigl\{ \tilde{\Psi}(\gamma) P_k \tilde{\Psi}^\top(\gamma) \bigr\}
\nonumber \\
&\quad + \beta^N X_{N|k} + \frac{\beta^{N+1}}{1-\beta} \EE\bigl\{ \tilde{D}(\gamma)
%\diag\{\Sigma_v,\Sigma_w\}
\Bigl[\begin{smallmatrix}\!\Sigma_v\! & \\ & \!\Sigma_w\! \end{smallmatrix}\Bigr]
\tilde{D}^\top(\gamma) \bigr\} .
\label{eq:term_cost_sdp}
\end{align}
Using (\ref{eq:xN_covar}) and Schur complements, (\ref{eq:term_cost_sdp}) can be expressed as a linear matrix inequality in $\theta_k = ({\bf c}_k, {\bf L}_k)$ and $P_k$.

Re-writing the constraints of problem (\ref{eq:opt}) using the matrix, $X_{i|k}$, of second moments yields the condition
\vspace{-2mm}
\[
\sum_{i=0}^\infty
\beta^i
\tr\bigl[ ( \ones_{2\times 2}
\otimes H^\top H)  X_{i|k}\bigr]
\leq \mu_{k},
\vspace{-1.5mm}
\]
which is equivalent to the constraint
\begin{equation}\label{eq:discounted_constraint}
\tr({\bf H}_\beta {\bf X}_k)
+
\tr\bigl[ ( \ones_{2\times 2} \otimes H^\top H)  P_k\bigr]
% \tr \Bigl[\begin{smallmatrix} \Xi^\top \Xi & \Xi^\top \Xi \\ \Xi^\top \Xi & \Xi^\top \Xi \end{smallmatrix}\Bigr]  P_k\Bigr)
% \Bigl[ \Bigl(\begin{bmatrix} I & I \\ I & I \end{bmatrix} \otimes \Xi^\top \Xi \Bigr)  P_k\Bigr]
\leq \mu_{k}
\end{equation}
where ${\bf H}_\beta = \ones_{2\times 2} \otimes \diag \{ H^\top H, \beta H^\top H, \ldots, \beta^{N-1} H^\top H \}$.
%where $\bs{\Xi} = \bigl[\begin{smallmatrix} 1 & 1 \\ 1 & 1 \end{smallmatrix}\bigr] \otimes \diag \{ \Xi^\top\Xi, \beta \Xi^\top\Xi, \ldots, \beta^{N-1} \Xi^\top\Xi \}$.
%

The expressions for the cost and constraints in (\ref{eq:cost_bound})-(\ref{eq:discounted_constraint}) allow the optimization~(\ref{eq:opt}) defining $\theta_k^\ast$ to be formulated as
%\begin{equation} \label{eq:opt2}
%\begin{aligned}
%\theta^\ast_k = \arg\min_{\theta_k,P_k}
%& \tr (\Q_\beta{\bf X}_k) \!+\! \tr(\R_\beta{\bf U}_k) \!+\!
%\tr\Bigl(\Bigl[\begin{smallmatrix} \!Q & Q\ \\ \!Q & Q + K^\top\! R K\! \end{smallmatrix}\Bigr]\!
%P_k\!\Bigr)
%\\
%\text{s.t.} \
%& \begin{aligned}[t]
%&
%\tr({\bf H}_\beta {\bf X}_k) +
%\tr \bigl[ (\ones_{2\times 2} \otimes H^\top H)  P_k\bigr] \leq \mu_{k} ,
%\\
%& P_k \succeq \beta \EE\bigl\{ \tilde{\Psi}(\gamma) P_k \tilde{\Psi}^\top(\gamma) \bigr\} + \beta^N X_{N|k} \\
%& \qquad+ \frac{\beta^{N+1}}{1-\beta} \EE\bigl\{ \tilde{D}(\gamma) \Bigl[\begin{smallmatrix} \!\Sigma_v\! & \\ & \!\Sigma_w\!\end{smallmatrix} \Bigr] \tilde{D}^{\top}\!(\gamma) \bigr\} .
%\end{aligned}
%\end{aligned} 
%\end{equation}
\begin{gather} 
\theta^\ast_k = \arg\min_{\theta_k,P_k}
\tr (\Q_\beta{\bf X}_k) \!+\! \tr(\R_\beta{\bf U}_k) \!+\!
\tr\Bigl(\Bigl[\begin{smallmatrix} \!Q & Q\ \\ \!Q & Q + K^\top\! R K\! \end{smallmatrix}\Bigr]\!
P_k\!\Bigr) \nonumber
\\
\text{s.t.} \  \eqref{eq:term_cost_sdp}, \eqref{eq:discounted_constraint}. \label{eq:opt2}
\end{gather}
% where ${\bf X}_k$, ${\bf U}_k$ and $X_{N|k}$ are defined by (\ref{eq:x_covar}), (\ref{eq:u_covar}) and (\ref{eq:xN_covar}).
\vspace{-5mm}
\begin{rem}
Problem~(\ref{eq:opt2}) can be expressed  as a semidefinite program in the variables $\theta_k=\{{\bf c}_k,{\bf L}_k\}$ and $P_k$ using (\ref{eq:2nd_pred_x_u}b), (\ref{eq:x_covar}) and (\ref{eq:xN_covar}).
Alternatively, we can eliminate $P_k$ from (\ref{eq:opt2}) by writing  the solution of the Lyapunov equation $P_k = \beta\EE\bigl\{\tilde{\Psi}(\gamma) P_k \tilde{\Psi}^\top(\gamma)\bigr\} + \Xi$ for given $\Xi=\Xi^\top$ as
\vspace{-2mm}
\[
\vvec(P_k) \!=\! \bigl[ I - \beta (1-\lambda) \tilde{\Psi}(0)\otimes\tilde{\Psi}(0) - \beta\lambda\tilde{\Psi}(1)\otimes\tilde{\Psi}(1)\bigr]^{-1}
\!\!\vvec(\Xi)
\vspace{-0.5mm}
\]
(where the matrix inverse can be computed of\mbox{}f\mbox{}line), which allows problem~(\ref{eq:opt2}) to be expressed, using standard matrix vectorization identities, as a convex quadratic program in $\theta_k = ({\bf c}_k,{\bf L}_k)$ with a single quadratic constraint.
%
% Alternatively, by using the fact that the Lyapunov equation $P_k = \beta\EE\bigl\{\tilde{\Psi}(\gamma) P_k \tilde{\Psi}^\top(\gamma)\bigr\} + \Xi$
% for any given matrix $\Xi=\Xi^\top$ has the explicit solution
% \[
% \vvec(P_k) = \bigl[ I - \beta (1-\lambda) \tilde{\Psi}(0)\otimes\tilde{\Psi}(0) - \beta\lambda\tilde{\Psi}(1)\otimes\tilde{\Psi}(1)\bigr]^{-1}
% \!\!\vvec(\Xi)
% \]
% where the matrix inverse can be computed offline, problem~(\ref{eq:opt2}) can be expressed as a convex quadratic program in $\theta_k = ({\bf c}_k,{\bf L}_k)$ with a single quadratic constraint.
\end{rem}

\section{Closed loop properties} \label{sec:properties}

This section considers the performance of the closed loop system~(\ref{eq:system}) with the control law of Section~\ref{sec:control_implementation}.
We use the solution $\theta^\ast_k=\{{\bf c}^\ast_k,{\bf L}^\ast_k\}$ of (\ref{eq:opt}) at time $k$ to construct a feasible, but possibly suboptimal, solution for (\ref{eq:opt}) at time \mbox{$k+1$} (i.e.~given $\U_k$, $\I_k$), which we denote  $\theta^\tail_{k+1}= \{{\bf c}^\tail_{k+1},{\bf L}^\tail_{k+1}\}$, where
\begingroup
\setlength{\abovedisplayskip}{-0.5pt}
\setlength{\belowdisplayskip}{0pt}
\begin{subequations}\label{eq:tail}
\begin{align}
{\bf c}^\tail_{k+1} &\defeq
\begin{bmatrix}
c^\ast_{1|k} \\ \vdots \\  c^\ast_{N-1|k} \\ 0
\end{bmatrix}
+
\begin{bmatrix}
L^\ast_{1,0|k} \\
 \vdots \\
L^\ast_{N-1,0|k} \\
0
\end{bmatrix} \gamma_k (y_k - C\xhat_k )  ,
\label{eq:c_tail}
\\
{\bf L}^\tail_{k+1} &\defeq
\begin{bmatrix}
L^\ast_{1,1|k} & & & \\
\vdots & \ddots & & \\
L^\ast_{N-1,1|k} & \cdots & L^\ast_{N-1,N-1|k} & \\
0 & \cdots & 0 & 0
\end{bmatrix} .
\end{align}
\end{subequations}
\endgroup
Following \cite{yan18}, we define the constraint threshold $\mu_{k}$ in~(\ref{eq:opt}) for all $k>0$ in terms of $\theta^\tail_{k}$. This ensures recursive feasibility of the MPC optimization without requiring bounds on the noise $v_k$ and disturbance $w_k$.
Thus
\begin{equation}\label{eq:mu_def}
\mu_k \defeq \begin{cases} \epsilon , & k = 0
\\
\tr ( {\bf H}_\beta {\bf X}^\tail_{k}) +
\tr \Bigl[ (\ones_{2\times 2} \otimes H^\top H)  P^\tail_k \Bigr] ,  & k > 0
\end{cases}
\end{equation}
where
%$\U_{k-1}$, $\I_{k-1}$ are assumed known when $\mu_k$ is computed and
\begin{align*}
P^\tail_k &= \beta \EE\bigl\{ \tilde{\Psi}(\gamma) P^\tail_k \tilde{\Psi}^\top(\gamma) \bigr\}
\\
& \quad + \beta^N X^\tail_{N|k} + \frac{\beta^{N+1}}{1-\beta} \EE\bigl\{ \tilde{D}(\gamma) \Bigl[\begin{smallmatrix} \Sigma_v & \\ & \Sigma_w\end{smallmatrix} \Bigr] \tilde{D}^{\top}(\gamma) \bigr\},
\\
{\bf X}^\tail_k &=
\begin{bmatrix} 0 & 0 \\ 0 & \bs{\pi}^\tail_k \bs{\pi}_k^{\tail\,\top} \end{bmatrix}
+
\begin{bmatrix} I & 0 \\ 0 & \bs{\Pi}_k^\tail\end{bmatrix} \bs{\Omega}_k \begin{bmatrix} I & 0 \\ 0 & \bs{\Pi}_k^\tail\end{bmatrix}^\top,
\\
X_{N|k}^\tail &=
\begin{bmatrix} 0 & 0 \\ 0 & \pi_{N|k}^\tail \pi_{N|k}^{\tail\,\top} \end{bmatrix}
+
\begin{bmatrix} I & 0 \\ 0 & \Pi_{N|k}^\tail\end{bmatrix} \Omega_{N|k} \begin{bmatrix} I & 0 \\ 0 & \Pi_{N|k}^\tail\end{bmatrix}^\top
\end{align*}
with
$\bs{\pi}_k^\tail = {\bf S}_{\Phi} \xhat_k + {\bf T}_{(\Phi,B)}{\bf c}_k^\tail$,
$\bs{\Pi}_k^\tail = {\bf T}_{(\Phi,B)} {\bf L}_k^\tail + {\bf T}_{(\Phi,A)}{\bf M}$.

% where
% ${\bf X}^\tail_{k}$ and $X^\tail_{N|k}$ are defined by (\ref{eq:x_covar}) and (\ref{eq:xN_covar}) with $\bs{\pi}_k$ and $\bs{\Pi}_k$ replaced by
% $\bs{\pi}_k^\tail = {\bf S}_{\Phi} \xhat_k + {\bf T}_{(\Phi,B)}{\bf c}_k^\tail$ and
% $\bs{\Pi}_k^\tail = {\bf T}_{(\Phi,B)} {\bf L}_k^\tail + {\bf T}_{(\Phi,A)}{\bf M}$.

\begin{thm}\label{thm:cl_constraint}
If problem~(\ref{eq:opt}) is feasible at $k=0$, then (\ref{eq:opt}) remains feasible for all $k>0$ and the state of (\ref{eq:system}) under the control law of Section~\ref{sec:control_implementation} satisfies
\vspace{-0.5mm}
\begin{equation}\label{eq:cl_constraint}
\sum_{k=0}^\infty \beta^k \EE \{ \| H x_{k}\|^2\} \leq \epsilon  .
\end{equation}
\end{thm}
\vspace{-2.5mm}
\bpf
The definition (\ref{eq:mu_def}) of $\mu_k$ trivially ensures feasibility for all $k>0$.
The definitions (\ref{eq:tail}a,b) ensure that,
%prior to receiving the sensor information $(\gamma_k,y_k)$,
at time $k$ (given $\U_{k-1}$, $\I_{k-1}$),
the distributions of the state and control sequences $\{x_{i|k+1}\}_{i=0}^\infty$ and $\{u_{i|k+1}\}_{i=0}^\infty$ are identical to the distributions of $\{x_{i+1|k}\}_{i=0}^\infty$ and $\{u_{i+1|k}\}_{i=0}^\infty$.
%considering $\v_k$, $w_k$ and $\gamma_k$ to be random variables,
%
%For any realisation of the estimation error $x_k-\xhat_k$ and disturbance and noise sequences ${\bf w}_k$ and ${\bf v}_k$, this choice of $\theta^\tail_{k+1}$ ensures $x_{i|k+1}=x_{i+1|k}$, $\xhat_{i|k+1} = \xhat_{i+1|k}$ and $u_{i|k+1} = u_{i+1|k}$, for all $i\geq 0$.
%
Therefore
\vspace{-3mm}
\[
\sum_{i=0}^\infty \beta^i \EE_k \{ \|H x_{i|k}\|^2 \}
=
\tr({\bf H}_\beta {\bf X}_k)
+
\tr \bigl[ (\ones_{2\times 2} \otimes H^\top H)  P_k \bigr]
\]
implies
\[
\beta \EE_k\{ \mu_{k+1} \}
\leq \mu_k - \EE_k\{ \|Hx_{0|k}\|^2 \}
= \mu_k - \EE_k\{ \|Hx_{k}\|^2 \} .
\]
Hence the trajectories of the closed loop system satisfy
\[
\sum_{i=0}^\infty \beta^{i} \EE_k \{ \|H x_{k+i}\|^2 \}
\leq
\mu_k - \lim_{i\to\infty}\beta^i \EE_{k} \{ \mu_{k+i} \}
\leq
\mu_k
\]
for all $k \geq 0$.
% and Chebyshev's inequality implies (\ref{eq:cl_constraint}).
\epf

\begin{cor}\label{cor:cl_constraint}
Let $J_k:=J(\theta_k^\ast,\xhat_k,\Sigma_k)$ denote the optimal value of the objective in (\ref{eq:opt}). Then under the control law of Section~\ref{sec:control_implementation}, the trajectories of (\ref{eq:system})  satisfy
\vspace{-1mm}
\begin{equation}\label{eq:cl_cost}
\sum_{k=0}^{\infty} \beta^k \EE \bigl\{ \| x_{k} \|_Q^2  + \| u_{k}\|_R^2 \bigr\} \leq J_0 .
\end{equation}
\end{cor}
\vspace{-2mm}
\bpf
Applying the same argument used in the proof of Theorem~\ref{thm:cl_constraint}
to the definition of the objective in (\ref{eq:opt}) yields
% \[
% J(\theta_k,\xhat_k,\Sigma_k) = \sum_{i=0}^{\infty} \beta^i \EE \bigl\{ \| x_{i|k} \|_Q^2 + \| u_{i|k}\|_R^2 \bigr\}
% \]
\begin{equation*}
\beta \EE_k \{ J(\theta^\tail_{k+1},\xhat_{k+1},\Sigma_{k+1})\} 
= J_k - \EE_k\{\|x_{k}\|_Q^2 + \|u_{k}\|_R^2\} ,
\end{equation*}
and since $ J_k \leq J(\theta^\tail_{k},\xhat_{k},\Sigma_{k}) ~\forall k$ by optimality, the bound in (\ref{eq:cl_cost}) follows.
\epf
%---

\section{Numerical examples}\label{sec:example}
This section gives a numerical example to demonstrate that the closed loop system satisfies \eqref{eq:cl_constraint} and \eqref{eq:cl_cost} and to compare with the unconstrained optimal LQG controller. We consider a system obtained by discretising a linearised continuous time model of a double inverted pendulum with a sample time of 0.01\,s as  in \citep{kwakernaak1985}.   The system matrices are
%\vspace{-1mm}
\begin{gather*}
A=
\left[\begin{smallmatrix}
    1.0005  &  0.01 &  -0.0005  &  0 \\
    0.098  &  1.0005 &  -0.0981  & -0.0005 \\
   -0.0005  &  0 &   1.0015  &  0.01 \\
   -0.0981  & -0.0005 &   0.2942  &  1.0015
\end{smallmatrix}\right], \\
B =
\left[\begin{smallmatrix}
    0.0001 &  -0.0001\\
    0.01 &  -0.02    \\
   -0.0001 &   0.0003\\
   -0.02 &   0.05
\end{smallmatrix}\right],
\quad
C=
\left[\begin{smallmatrix}
     1 &    0 &    0 &    0\\
     0 &    0 &    1 &    0
\end{smallmatrix}\right],
\quad
D=I,
\end{gather*}%
and $\omega_k \sim \mathcal{N}(0,\Sigma_w)$, $v_k \sim \mathcal{N}(0,\Sigma_v)$, $\lambda=0.6$. Here $\Sigma_w=\diag\{0.5, 0.2 , 0.9, 0.3\}$ and $\Sigma_v=1.1 I$. Initial conditions are given by % $x_0=[-0.8, 0.4, 0.55, -0.5]^\top$, $\xhat_0=[0.1, 0.05, 0.1, 0.05]^\top$ and
%\vspace{-1mm}
%{\small
\begin{equation*}
x_0=\left[\begin{smallmatrix}  -0.8\\ 0.4\\ 0.55\\ -0.5 \end{smallmatrix}\right], \ 
\xhat_0= \left[\begin{smallmatrix} 0.1\\ 0.05\\ 0.1\\ 0.05 \end{smallmatrix}\right], \ 
\Sigma_0=
\left[\begin{smallmatrix}
    0.5 &  -0.5 &  -0.5 &  0.5\\
   -0.5 &   0.5 &  0.5  & -0.5\\
   -0.5 &   0.5 &  0.5  & -0.5\\
    0.5 &  -0.5 & -0.5  &  0.5
\end{smallmatrix}\right].
\end{equation*}%}%
The constraint of \eqref{eq:orig_opt_problem} is defined by $\beta=0.95$, $\epsilon=111$ and
% \begingroup
% \setlength{\abovedisplayskip}{-0.5pt}
% \setlength{\belowdisplayskip}{-0.5pt}
%{\small
\[H=
\left[ \begin{smallmatrix}
0  &  0.1  &  0 &  -0.1\\
0.1  &  0  & -0.1 &  0
\end{smallmatrix} \right].
\]% }%
% \endgroup
The weighting matrices in the cost function of \eqref{eq:orig_opt_problem} are given by
$Q=\diag \{10, 0.1, 10, 0.1\}$, $R=10^{-4} I$. We choose a prediction horizon $N=5$, $K$ as the unconstrained LQ-optimal, $K_{LQ}$, with respect to $(A,B,Q,R)$ and $M=\bar{\Sigma}C^\top (C \bar{\Sigma}C^\top + \Sigma_v)^{-1}$, where $\bar{\Sigma}$ is the solution of the algebraic Riccati equation
\[
\bar{\Sigma}=A \bar{\Sigma} A^\top +\Sigma_w -\lambda A\bar{\Sigma}C^\top (C\bar{\Sigma}C^\top +\Sigma_v)^{-1}C\bar{\Sigma}A^\top.
\]
Using the above information, we solve problem \eqref{eq:opt2} and obtain $J_0=2.368\times 10^4$.

\textit{Simulation A}: To estimate empirically the LHS of \eqref{eq:cl_constraint} and \eqref{eq:cl_cost}, we
consider their average values over $10^3$ simulations, each of which has a length of $500$ time steps. This gives $\sum_{k=0}^{\infty} \beta^k \EE \{ \| H x_{k}\|^2\}$ and $\sum_{k=0}^{\infty} \beta^k \EE \bigl\{ \| x_{k} \|_Q^2  + \| u_{k}\|_R^2 \bigr\}$ as $104.7$ and $4.774\times 10^3$ respectively. Therefore, these estimates agree with the bound \eqref{eq:cl_constraint} and \eqref{eq:cl_cost}. Moreover,
$\beta^{500} = 7.3 \times 10^{-12}$, so a further increase in the horizon length has negligible
effect on these estimates.

\textit{Simulation B}: To compare with the above results, we run the same number of simulations with the same $\{\omega_k\}$, $\{v_k\}$,$\{\gamma_k\}$ sequences using the unconstrained optimal LQG controller, where $u_k=K_{LQ}\xhat_k$ and the estimator gain is time-varying and given by $M=\Sigma_k C^\top (C \Sigma_k C^\top + \Sigma_v)^{-1}$. Here $\Sigma_k$ evolves as
\[
\Sigma_{k+1}\!=\!A \Sigma_k A^\top +\Sigma_w -\gamma_k A\Sigma_k C^\top (C\Sigma_k C^\top +\Sigma_v)^{-1}C\Sigma_kA^\top.
\]
This gives $\sum_{k=0}^{\infty} \beta^k \EE \{ \| H x_{k}\|^2\}$ as $123.8$, violating the bound \eqref{eq:cl_constraint}, and
$\sum_{k=0}^{\infty} \beta^k \EE \bigl\{ \| x_{k} \|_Q^2  + \| u_{k}\|_R^2 \bigr\}$ as $3.626 \times 10^3$, which is smaller than that in Simulation A, as expected.

\section{Conclusion}\label{sec:conclusion}
This paper describes an output feedback MPC algorithm for linear discrete time systems with additive disturbances and noisy sensor measurements transmitted over a packet-dropping communication channel.
%By employing a discounting factor, we do not require that disturbances and measurement noise be bounded. 
%
By designing a control policy with an affine dependence on future observations, we provide a convex formulation of a stochastic quadratic regulation problem subject to a discounted expectation constraint.
Our controller parameterization ensures recursive feasibility of the MPC optimization problem and ensures a cost bound and constraint satisfaction in closed loop operation.
Future work will explore interconnections between conditions for mean square stability of the MPC law and the values of the packet loss probability and the discount factor in the receding horizon optimization.

\begin{ack}
The authors would like to thank Prof. Subhrakanti Dey for helpful discussions about early versions of this paper.
\end{ack}

\balance 
\bibliography{ifacconf}             % bib file to produce the bibliography
                                                  % with bibtex (preferred)

% \appendix
% \section{A summary of Latin grammar}    % Each appendix must have a short title.
% \section{Some Latin vocabulary}              % Sections and subsections are supported
%                                                                          % in the appendices.

\end{document}